\documentclass{amsart}
\usepackage[dvips]{epsfig}

\newcommand{\B}{\mathcal B}
\newcommand{\C}{\mathbb C}
\newcommand{\CP}{\mathbb{CP}^1}
\newcommand{\dpqr}{\Delta(p,q,r)}
\newcommand{\fl}[1]{\lfloor #1 \rfloor}
\newcommand{\fltwo}[1]{\lfloor \frac{#1}{2} \rfloor}
\newcommand{\pM}{\partial M}
\newcommand{\pt}{\tilde{\pi}}
\newcommand{\PSLC}{\mbox{PSL}_2(\C)}
\newcommand{\Q}{\mathbb Q}
\newcommand{\R}{\mathbb R}
\newcommand{\SLC}{\mbox{SL}_2(\C)}
\newcommand{\St}{\Sigma_2}
\newcommand{\tXi}{\widetilde{X}_i}
\newcommand{\Z}{\mathbb Z}

\begin{document}
\title{The Culler-Shalen seminorms of the $(-3,3,4)$ pretzel knot}

\author{Thomas W.\ Mattman}
\address{Department of Mathematics, McGill University\\
Montr\'eal, Qu\'ebec, Canada}
\email{mattman@math.mcgill.ca}
\thanks{Research supported by grants NSERC OGP 0009446 and FCAR EQ 3518.}

\dedicatory{Dedicated to Professor Murasugi on the occasion of his
70th birthday.}

\begin{abstract}
We describe a method to compute the Culler-Shalen seminorms
of a knot, using the $(-3,3,4)$ pretzel knot as an illustrative
example. We deduce that the $\SLC$-character variety of this 
knot consists of three algebraic curves and that it admits
no non-trivial cyclic or finite surgeries. We also summarize
similar results for other $(-3,3,n)$ pretzel knots
including the observation that the Seifert surgeries for these
knots are precisely those integral slopes lying between two
of the boundary slopes.
\end{abstract}

\maketitle

\section*{Introduction}
\begin{figure}
\begin{center}
\epsfig{file=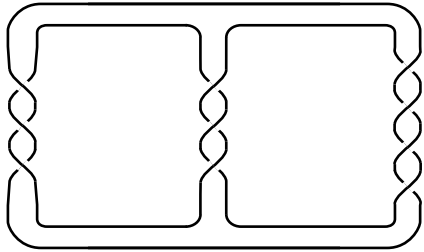}
\end{center}
\caption{The $(-3,3,4)$ pretzel knot \label{fg334}}
\end{figure}

Let $K \subset S^3$ be a hyperbolic knot. We have been developing
machinery which parlays a little information about such a knot into
a rather thorough understanding of its $\SLC$ character variety as well
as a classification of its finite and cyclic fillings. As an example
of our method, we will examine the case where $K$ is the  
$(-3,3,4)$ pretzel knot (see Figure~\ref{fg334}).

There are three types of input necessary for our approach.
\begin{enumerate}
\item A listing of the boundary slopes of the knot. For example
Hatcher and Oertel~\cite{HO} have discussed how to determine
the boundary slopes of any Montesinos knot.
\item Information about $\St$, the two-fold branched cyclic
cover of the knot. Again, Montesinos knots are good candidates
in this regard as $\St$ is then a Seifert fibred manifold.
\item At least one non-trivial finite, cyclic, or small Seifert surgery.
\end{enumerate}
Given this data, we can generally determine the number
of components of the $\SLC$ character variety of the knot as
well as the Culler-Shalen seminorms on each component. This
allows us to classify all finite and cyclic surgeries on the knot.
Part of our motivation in presenting this work is to solicit
examples of other knots meeting our input criteria.

Our algorithm has been successfully applied to 
the twist knots~\cite{BMZ} and the $(-2,3,n)$ pretzel 
knots~\cite{M1}. The twist knots are Montesinos knots which each admit 
three small Seifert surgeries and are therefore amenable to our methods.
The $\SLC$ character variety of a twist knot consists of two 
algebraic curves and these knots admit no non-trivial finite or
cyclic surgeries. As for the $(-2,3,n)$ pretzel knots, since 
they each have
two small Seifert surgeries, we can use our techniques to show
that the character variety consists of two or three curves and
that there are only five non-trivial finite surgeries. That is,
the $(-2,3,7)$ pretzel has three non-trivial finite surgeries, the $(-2,3,9)$
has two non-trivial finite surgeries and the remaining hyperbolic $(-2,3,n)$
pretzel knots admit no non-trivial finite or cyclic surgeries. 

To illustrate our methods, we will look at the $(-3,3,4)$ pretzel
knot which has a small Seifert surgery of slope 1. Our intention is
to give an overview of the main ideas of our approach. For a more
careful account we refer the reader to \cite{BMZ,M1,M2}.

\section*{The character variety and Culler-Shalen seminorms}

Our main tool is the Culler-Shalen seminorm which we now briefly
describe. A more detailed exposition can be found in \cite[Chapter 1]{CGLS}
or \cite{BZ2}.

Let $R = \mbox{Hom}(\pi, \SLC)$ denote the set of
$\SLC$-representations of the fundamental group $\pi$ of 
$M = S^3 \setminus K$.
Then $R$ is an affine algebraic set, as is $X$, the set of characters
of representations in $R$.
Since $M$ is small~\cite{O}, the irreducible components of $X$ are 
curves~\cite[Proposition 2.4]{CCGLS}. Moreover, for each component
$R_i$ of $R$ which contains an irreducible representation, the 
corresponding curve $X_i$ induces a non-zero seminorm $\| \cdot \|_i$ on 
$V = H_1(\pM; \R)$~\cite[Propositon 5.7]{BZ2} via the following construction.

For $\gamma \in \pi$, define the regular function
$I_{\gamma}:X \to \C$ by $I_{\gamma}(\chi_{\rho}) = \chi_{\rho}(\gamma) =
\mbox{trace}(\rho(\gamma))$.
By the Hurewicz isomorphism,
a class $\gamma \in L =  H_1(\pM ;{\Z})$ determines an element
of $\pi_1(\pM)$, and therefore an element of $\pi$ well-defined
up to conjugacy. The function
$f_{\gamma} = I_{\gamma}^2 -4$ is again regular and so can be pulled back to
$\tXi$, the smooth projective variety birationally equivalent to
$X_i$. For $\gamma \in L$, $\| \gamma \|_i$ is the degree of
$f_{\gamma} : \tXi \to \CP.$ The seminorm is extended to $V$ 
by linearity. We will call a seminorm constructed in this manner 
a Culler-Shalen seminorm.

If no $f_{\gamma}$ is constant on $X_i$, then $\| \cdot \|_i$ is in fact a
norm and we shall refer to $X_i$ as a {\it norm curve}. If $X_i$ is not 
a norm curve, then there is a boundary slope 
$r$ such that $f_r$ is constant on $X_i$.
In this case, we will call $X_i$ an {\it $r$-curve}.
The minimal norm 
$s_i = \min \{ \| \gamma \|_i \, ; \, \gamma \in L, \, \| \gamma \|_i > 0 \}$ 
is an even integer, as is
$S = \sum_{i} s_i$, the sum being taken over the curves $X_i \subset X$. 

Our goal is to show that the character variety of the $(-3,3,4)$ 
pretzel knot has one norm curve and one $r$-curve with $r=0$.
As the reducible characters also form a curve, this means $X$
consists of exactly three algebraic curves. We will present
the argument by discussing how each of the three inputs mentioned 
in the introduction come into play.

\section*{The two-fold branched cyclic cover}

As the $(-3,3,4)$ pretzel knot is a Montesinos knot, 
$\St$ is a Seifert fibred space. The base orbifold $\B$ of $\St$ is
$S^2$ with cone points of order $3,3$ and $4$. We can use 
this information to determine $S$. Essentially, the 
argument relies on the strong connections between
the various fundamental groups. To wit, let $\pt$
be the index two subgroup of $\pi$ corresponding to 
the two-fold cyclic cover. Then $\pi_1(\St) = \pt/$ $<\mu^2>$
where $\mu$ is the class of a meridian of $K$. Also, 
$\pi_1^{\mbox{orb}}(\B) = \Delta(3,3,4) = $ $<a,b | a^3, b^3, (ab)^4>$
is a triangle group and isomorphic to $\pi_1(\St)/Z(\pi_1(\St))$,
where $Z(\cdot)$ denotes the center.

Now, by~\cite[Corollary 1.1.4]{CGLS}, $\| \mu \|_i = s_i$ for
each curve $X_i \subset X$. So
\begin{eqnarray*}
S & = & \sum \| \mu \|_i \\
  & = & \sum \| 2 \mu \|_i - \| \mu \| \\
  & = & \sum Z_x(f_{\mu^2}) - Z_x(f_{\mu}),
\end{eqnarray*}
where $Z_x(\cdot)$ denotes the degree of zero of the function
at the point $x \in X$.
To evaluate the sum, we are led to look at zeroes of $f_{\mu^2}$
which are not zeroes of $f_{\mu}$. If $0 = f_{\mu^2} = I_{\mu^2} - 4$,
then $\mbox{trace}(\rho(\mu^2)) = \pm 2$. In other words, the
representations which kill $\mu^2$, and therefore factor through
$\pi_1(\St)$, will all contribute to the sum. 

The abelian $\SLC$-representations of $\pi_1(\St)$ lift to
(binary) dihedral representations of $\pi$. The number of 
$\SLC$-characters of these dihedral representations may
be related to the Alexander polynomial $\Delta_K(t)$ \cite[Theorem 10]{K}
and there are 
($(|\Delta_K(-1)| - 1)/2 =$) 4 such $\SLC$ characters. 

On the other hand, non-abelian representations will factor through
the center of $\pi_1(\St)$ to become representations of
$\Delta(3,3,4)$. In general, 
the number of $\PSLC$-characters of $\dpqr$ is (see \cite[Proposition 3.2]{BB})
\begin{eqnarray} \label{eqpqr}
&& (p- \fltwo{p} - 1)(q- \fltwo{q} - 1)(r- \fltwo{r} - 1) +
   \fltwo{p} \fltwo{q} \fltwo{r} \\
&& \mbox{}+ \fltwo{\gcd(p,q)} + \fltwo{\gcd(p,r)} + \fltwo{\gcd(q,r)}
   + 1 \nonumber \end{eqnarray}
where $\fl{x}$ denotes the largest integer less than or equal to $x$.
This count includes the reducible characters. As the character of a 
reducible representation is also the character of an abelian representation,
we see that the reducible characters correspond to representations of
$H_1(\dpqr) = \Z / a \oplus \Z / (b/a)$ where $a = \gcd(p,q,r)$
and $b = \gcd(pq,pr,qr)$. So the number of reducible 
$\PSLC$-characters of $\dpqr$ is 
\begin{equation} \label{eqH1}
   \begin{array}{cl} \fltwo{b} + 1, & \mbox{if } a \equiv 1 \pmod{2}  \\
                     \fltwo{b} + 2, & \mbox{if } a \equiv 0 \pmod{2}.
   \end{array}
\end{equation}

In particular, $\Delta(3,3,4)$ admits $3$ irreducible $\PSLC$
characters
and therefore $6$ irreducible $\SLC$ characters \cite[Lemma 5.5]{BZ1}.

By \cite[Theorem A]{BB},
each of the dihedral and triangle group characters contributes two 
to $S$ so that $S = 2(6 + 4) = 20$.

\section*{Boundary slopes}

According to \cite{HO}, the boundary slopes of $K$ are $-14$, $0$ and $8/5$.
A small rearrangement of \cite[Lemma 6.2]{BZ1} shows that
\begin{equation} \label{eqgi}
\| \gamma \|_i = 2 \sum_j a_j^i \Delta(\gamma, \beta_j)
\end{equation}
the sum being taken over the boundary slopes $\beta_j$. Here, 
$\Delta(\gamma, \beta)$ denotes the minimal geometric intersection 
of curves representing $\gamma$ and $\beta$ in $\pi_1(\pM)$.
In particular, using standard meridian-longitude coordinates,
we can denote $\gamma$ (respectively $\beta$) as $a/b$ ($c/d$) 
$\in \Q \cup \{ 1/0 \}$.
Then $\Delta( \gamma, \beta) = |ad - bc|$.

Thus, given a list of boundary slopes, finding the Culler-Shalen
seminorm comes down to solving a system of equations for the
non-negative integers $a_j^i$. To fix the ordering in what 
follows, let $\beta_1 = -14$, $\beta_2 = 0$ and $\beta_3 = 8/5$. 

\section*{Cyclic, finite, or small Seifert fillings}

In order to solve Equation~\ref{eqgi}, we will need to know
about at least one non-trivial cyclic, finite, or small Seifert
filling $\alpha$, as $\| \alpha \|_i$ may then be related
to $s_i$.
\begin{itemize}
\item If $\alpha$ is a cyclic surgery, then $\| \alpha \|_i = s_i$
 \cite[Corollary 1.1.4]{CGLS}.
\item If $\alpha$ is a finite surgery, then $\| \alpha \|_i \leq
\max(2 s_i, s_i + 8)$ \cite[Theorem 2.3]{BZ1}.
\item If $\alpha$ is small Seifert then $M(\alpha)$, the filling
along $\alpha$, is Seifert fibred. The base orbifold will
be $S^2$ with cone points of order $p,q$ and $r$. In this case,
 $\| \alpha \|_i \leq s_i + C_{p,q,r}$
where $C_{p.q.r}$ is a constant depending on $p,q$ and $r$
(an example of this type follows).
\end{itemize}

For the $(-3,3,4)$ pretzel knot, $1$ filling is Seifert fibred with
base orbifold $S^2(2,5,7)$. Then the $6$ irreducible
$\PSLC$-characters of $\Delta(2,5,7)$ (see Equations~\ref{eqpqr} and
\ref{eqH1}) become $12$ $\SLC$-characters \cite[Lemma 5.5]{BZ1} 
each contributing $2$ to $\|1\|_i$ \cite[Theorem A]{BB}.
Thus $\sum \| 1 \|_i = S + 24$ (summing over the curves $X_i \subset X$).
Consequently, $s_i \leq \| 1 \|_i \leq s_i + 24$. 
Of course
$\| \mu \|_i = s_i$ since $M(\mu) = S^3$ is a cyclic filling. 

So we
have the equations
\begin{eqnarray}
 \| \mu \| = 2(a_1 + a_2 + 5 a_3)  & = & s \leq 20; \mbox{ and }
 \label{eq331}\\
 s \leq \| 1 \| = 2(15a_1 + a_2 + 3 a_3) & \leq & s + 24 \label{eq332}
\end{eqnarray}
(where we've suppressed the `$i$' sub- and superscripts).
Subtracting, we find 
\begin{equation}
0 \leq 7 a_1 - a_3 \leq 6. \label{eq333}
\end{equation}

Let us first investigate the case where $\| \cdot \|$ is 
a norm curve.
In order to have a norm (rather than just a seminorm), at least
two of the $a_j$ must be non-zero.
Equation~\ref{eq331} shows that $a_3 \leq 1$ on a norm curve. On
the other hand, from Equation~\ref{eq333} we see that $a_3 = 0$ 
implies $a_1 =0$ which would not be possible for a norm curve. 
Therefore
$a_3 =1$. Then Equation~\ref{eq333} implies $a_1 = 1$. Finally,
Equation~\ref{eq331} allows us to bound $a_2$: $0 \leq a_2 \leq 4$.
In particular, on a norm curve, we have $\| 1 \| = s + 24$. 
This means there can be
at most one norm curve. (If there were two or more, 
$\sum_i \| 1 \|_i  \geq S + 2(24)$ contradicting an earlier equation.)
On the other hand, the component $X_0$ of the character variety
containing the character corresponding to the holonomy representation 
{\it is} a norm curve \cite[Chapter 1]{CGLS}.
Let $\| \cdot \|_0$ denote the norm on $X_0$.

Since $\sum_i \| 1 \|_i = S + 24$, we see that $\| 1 \|_i = s_i$
for any $r$-curves. On the other hand, on an 
$r$-curve $X_i$, $\| 1 \|_i = s_i \Delta(1,r)$
and $r$ is a boundary slope \cite[Proposition 5.4]{BZ2}. 
So the only candidate is $r=0$. 

Now $M(0) = M_1 \cup M_2$, with 
$M_1$ Seifert fibred over $D^2 (2,2)$ and $M_2$ Seifert over 
$D^2(3,3)$, is a graph manifold 
and its $\PSLC$ representations will factor through
$\Z/2 \ast \Z/3$. Since the $\PSLC$-character variety
$\bar{X}(\Z/2 \ast \Z/3)$ contains exactly one curve 
\cite[Example 3.2]{BZ2}, the same is true of $\bar{X}(M(0))$ and
we conclude that there is a unique $r$-curve $X_1$ with $r=0$. 
Moreover the minimal norm is $s_1 = 2$. (For a more detailed 
account of this argument, see the discussion of the 
$M(2n+6)$ filling of the $(-2,3,n)$ pretzel knot in \cite{M1,M2}.)

\begin{figure}
\begin{center}
\epsfig{file=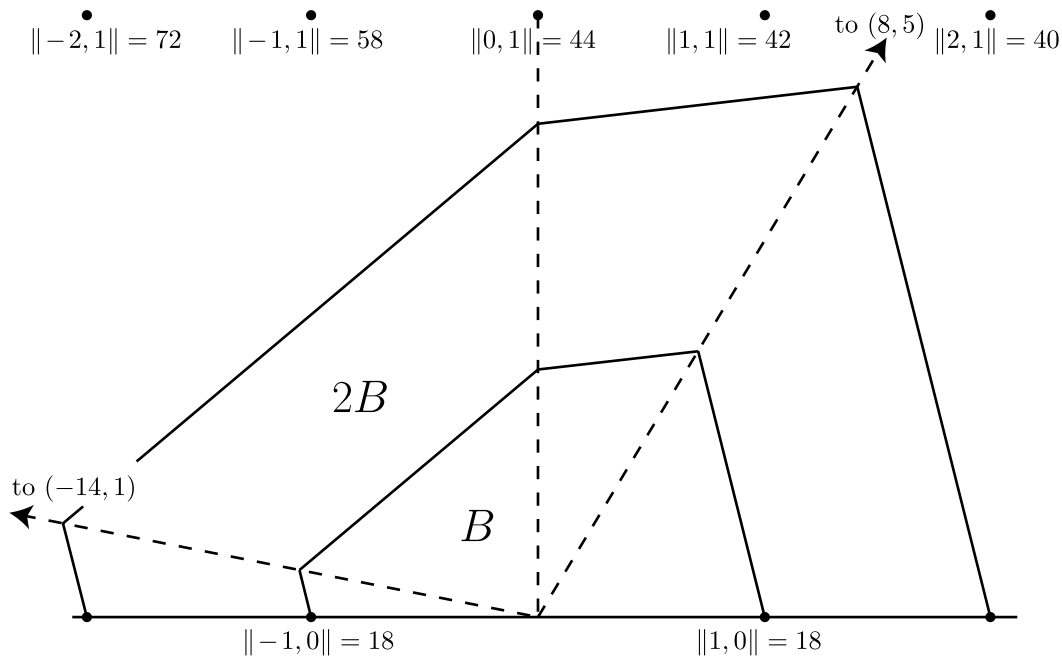}
\end{center}
\caption{The fundamental polygon $B$ of the $(-3,3,4)$ pretzel knot \label{fgB}}
\end{figure}

\begin{figure}
\begin{center}
\epsfig{file=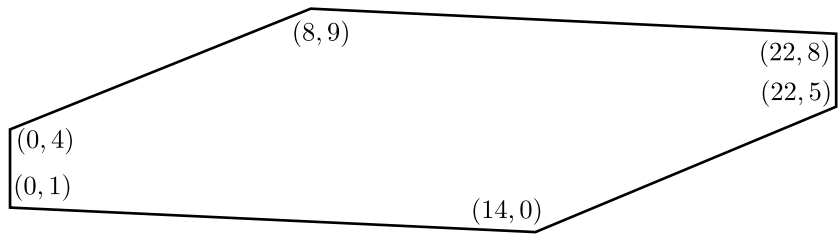}
\end{center}
\caption{The Newton polygon of the $A$ polynomial for the 
$(-3,3,4)$ pretzel knot \label{fgN}}
\end{figure}

Thus $X$ contains one norm curve $X_0$ and one $r$-curve $X_1$. The
Culler-Shalen seminorm on $X_1$ is $\| \gamma \|_1 = 2 \Delta(\gamma, 0)$
while that of $X_0$ is 
$$ \| \gamma \|_0 = 2 {[}\Delta(\gamma, -14) + 3 \Delta(\gamma, 0) + 
\Delta(\gamma, 8/5){]} $$
(i.e.\ we choose $a_2^0 = 3$ so that $\| \mu \|_0 + \| \mu \|_1 = S = 20$).
The polygon $B$ of radius $s_0 = 18$ in $\| \cdot \|_0$ is illustrated
in Figure~\ref{fgB}.
Notice that $B$ lies below the line $y=1/2$. Since any finite or
cyclic surgeries would have norm at most $\max(2s_0, s_0 +8)$
\cite[Theorem 2.3]{BZ1}, 
we see that $K$ admits no other cyclic or finite surgeries beyond 
trivial surgery along the meridian $\mu = 1/0$. (Delman~\cite{D}
has already shown that this knot has no non-trivial finite surgeries using 
completely different methods.)
The Newton polygon of
the $A$ polynomial is dual to $B$ \cite{BZ3} and is illustrated in
Figure~\ref{fgN}.

\subsection*{$(-3,3,n)$ pretzel knots} 

We now generalize to
the $(-3,3,n)$ pretzel knot which we will denote
by $K_n$. Note that $K_{-n}$ is the mirror reflection of $K_n$,
so we can assume $n \geq 0$. This family includes some knots
we have investigated elsewhere: $K_1$ is a twist knot \cite{BMZ} and $K_2$
is the reflection of the $(-2,3,-3)$ pretzel knot \cite{M1}. Since
$K_0$ is not prime, it's not hyperbolic 
and therefore not amenable to the techniques we have been discussing.
On the other hand, when $3 \leq n \leq 6$,
$K_n$ is hyperbolic. Moreover these knots have a 
Seifert surgery at slope $r =1$. (We have verified
this for $n=3,4,6$ using the Montesinos trick. For $n=5$ we have
only the evidence of SNAPPEA~\cite{Wk}.) So for these knots we have
the required inputs in order to 
apply our machinery and work out the Culler-Shalen seminorms.
 
But what of $n\geq 7$? Why stop at $n=6$? We are obliged to stop
since we have no evidence of $K_n$ admitting a Seifert filling
for $n \geq 7$. Indeed, the Seifert surgeries occur according to
a very nice pattern. 
By \cite{HO}, the boundary slopes of $K_n$ are $-(2n+6)$, $0$ and $8/(n+1)$. 
For $1 \leq n \leq 6$, the Seifert surgeries 
lie between the boundary slopes $0$ and $8/(n+1)$ as the 
following table illustrates.

\vspace{12pt}

\begin{center}
\begin{tabular}{|c|c|c|c|c|c|c|c|} \hline
$n$ & $1$ & $2$ & $3$ & $4$ & $5$ & $6$ & $\geq 7$\\ \hline
$8/(n+1)$ & $4$ & $8/3$ & $2$ & $8/5$ & $4/3$ & $8/7$ & $\leq 1$ \\ \hline
Seifert   &         &       &     &     &     &     &     \\
Surgeries & $1,2,3$ & $1,2$ & $1$ & $1$ & $1$ & $1$ & none \\ \hline
\end{tabular}
\end{center}

\vspace{12pt}

As the boundary slope $8/(n+1)$ moves across the integers toward
$0$, those integers cease to be available for Seifert surgeries.
For example, when $n \geq 7$, the boundary slope is $\leq 1$ and
there are no more Seifert surgeries. I should emphasize that this is
based on experimental evidence. These knots may admit other Seifert
surgeries beyond those I've listed in the table. In addition, although
I (or others) have shown that all the other surgeries in the table are
Seifert, the only evidence I have in the $n=5$ case comes from 
SNAPPEA~\cite{Wk}. Nonetheless, it is a curious pattern and
it would be nice to understand this phenomenon.

Thus we can only hope to apply our machinery to $K_n$ when $1 \leq n
\leq 6$. The first two cases are treated elsewhere and $n=4$ was 
discussed in detail above. 
For $K_3$ our method breaks down as
the equations corresponding
to Equations~\ref{eq331}, \ref{eq332} and \ref{eq333} above don't
result in a unique solution for the $a_j^i$'s. 
Since $K_5$ is not strongly invertible, we can not use the 
Montesinos trick to work out the indices for the Seifert surgery
of slope $1$. Without that information, we can not complete
the analysis of that knot.

However, $K_6$ is tractable. For this knot we
have the same conclusions as for $K_4$: there's one norm curve
with $s_0 = 22$ and 
$$\| \gamma \|_0 = 2 {[} \Delta( \gamma, -18) + 3 \Delta( \gamma, 0) + 
\Delta( \gamma, 8/7) {]},$$ 
and one $r$-curve with $r=0$ and $s_1 = 2$. This means that $K_6$ also
admits no non-trivial cyclic or finite surgeries. (Again,
Delman~\cite{D} had shown this previously using different methods.)

\section*{Acknowledgments}

I would like to thank Steve Boyer for many useful conversations and
in particular for suggesting that other pretzel knots might have 
Seifert fillings. I thank Professor Sakuma and the other organizers
for putting together a very nice workshop and for giving me the opportunity
to talk about my research.

\end{document}